\documentclass{llncs}
\usepackage{llncsdoc,color,cite,graphicx}
\usepackage{ActiveFlowDefn}

\begin{document}
\title{Model Reduction for DAEs with an Application to Flow Control
}

\titlerunning{Model Reduction with Application to Flow Control}

\author{J.T.~Borggaard  \and S.~Gugercin}
\institute{Department of Mathematics, Virginia Tech   \\
Interdiciplinary Center for Applied Mathematics, Virginia Tech\\
Blacksburg, VA, 24061-0123, USA  \\
\email{\{jborggaard,gugercin\}@vt.edu}
}

\maketitle

\begin{abstract}
Direct numerical simulation of dynamical systems is of fundamental importance
in studying a wide range of complex physical phenomena. However, the ever-increasing need for accuracy leads to extremely large-scale dynamical systems whose simulations impose huge computational demands. Model reduction offers one remedy to this problem by producing simpler reduced models that are both easier to analyze and faster to simulate while accurately replicating the original behavior.  Interpolatory model reduction methods have emerged as effective candidates for very large-scale problems due to their ability to produce high-fidelity (optimal in some cases) reduced models for linear and bilinear dynamical systems with modest computational cost. In this paper, we will briefly review the interpolation framework for model reduction and describe a well studied flow control problem that requires model reduction of a large scale system of differential algebraic equations.  We show that interpolatory model reduction produces a feedback control strategy that matches the structure of much more expensive control design methodologies.
\end{abstract}
\section{Introduction}
Direct numerical simulation of dynamical systems has been one of the few available means for studying many complex systems of scientific
and industrial value as dynamical systems are the basic framework for modeling, optimization, and control of
these complex systems. Examples include chemically reacting flows, fluid dynamics, and signal propagation and
interference in electric circuits. However, the ever increasing complexity and need for improved
accuracy lead to the inclusion of greater detail in the model, and inevitably finer discretizations.  Combined with the potential coupling to other complex systems, this results in extremely large-scale dynamical systems with millions of degrees of freedom to simulate. The simulations in these settings can be overwhelming; which is the main motivation for model reduction. The goal is to construct reduced models with significantly lower number of degrees of freedom that are easier to analyze and faster to simulate yet accurately approximate the important features in the underlying full-order large-scale 
simulations. 

There is a tremendous amount of literature on model reduction. Here we only include a partial list of various applications settings where model reduction has been applied with great success: In fluid flow \cite{sirovich1,Deane1991,Gatski1992,Ly2001,Hay09,Mathelin2009} and design of feedback control systems
\cite{AtwKin01,BBKT02,GABGcdc04,King98,Noack2003jfm,Noack2004acc,RCM00}, in optimization 
\cite{Arian2002,Kunisch2008,Antil2011,Antil2012,Baur11,Druskin2011solution,yue2013}, in nonlinear inverse problems  
 \cite{Wang:2005b,Galbally2009,Lieberman2010,Druskin2011solution,morpals2013,BuiThanh2008_AIAA}, in optimal design 
 \cite{Lieu2006,Lieu2007,Amsallem2008},  in the analysis of 
 structural mechanics \cite{su1991model,gawronski1990model,meyer1996balancing,bindel2006model,rudnyi2002review},
 in circuit theory \cite{FelF95,odabasioglu1997prima,Daniel2004,BonD07,Baur_etal2011,Fen05,bai1998make,reis2010positive,BenHtM11},
and in structural mechanics, such as \cite{su1991model,gawronski1990model,meyer1996balancing,bindel2006model,rudnyi2002review}. For a detailed discussion of several model reduction topics, see \cite{holmes,Ant05,BenHtM11,benner2005dimension,obinata2001model}.

The active flow control application we consider in this paper is a well studied flow control problem of stabilizing the von K\'arm\'an vortex shedding behind a circular cylinder by controlling the rotational velocity of the cylinder.  Upon linearizing the Navier-Stokes equations about a desired steady-state solution, the resulting large-scale linear systems of differential algebraic equations (DAE) is reduced by the interpolatory model reduction  framework recently developed in \cite{GSW2013}.  We use this reduced model to design the feedback control strategy and compare these results to other feedback control laws found in the literature.

\section{Description of the Flow Control Problem\label{sec:flow}} 
Suppression of the vortex shedding behind a bluff body is a classical
flow control problem with numerous applications ranging from minimizing
drag to reducing the cross-stream lift-induced fatigue cycling.  A number of experimental and computational studies have shown the effectiveness of cylinder surface suction and zero-mass actuation to completely suppress the von K\'arm\'an vortex street at modest Reynolds numbers that are slightly higher than the critical Reynolds number (bifurcation parameter) $Re_c\approx 47$, cf. \cite{park1994fcv,gunzburger1996fck}, based on the cylinder diameter and inflow velocity.

An alternate strategy capitalizes on the Magnus effect produced by cylinder rotation \cite{prandtl1925mew}.  In a sequence of experiments in \cite{taneda1978vof},
the authors showed that a rotationally oscillating cylinder using carefully selected choices of frequency, amplitude, and phase angle could {\em effectively eliminate the wake} for moderate flows of $Re\approx 61$ and $Re\approx 110$.  At higher values of the Reynolds number, it was not possible to eliminate the wake, but using good choices of frequency and amplitude made it possible to achieve nearly 20\% drag reduction.  Other experimental studies \cite{olinger1988ndw,filler1991rsl,tokumaru1991roc} suggest that matching the oscillation frequency to the vortex shedding frequency maximizes the impact on the flow (at this range of $Re$).  This was confirmed numerically in \cite{ou1991cof}.

A number of active feedback control approaches for the rotating/oscillating cylinder have appeared in the literature over the past fifteen years, including
\cite{he2000acd,xiao2002rcc,protas2004lfs,xiao2006ric,bewley2007mec,borggaard2008ros,stoyanov2009rom,akhtar2010crc,pralits2010isf}.

In the remainder of this section we describe our feedback control
strategy based on linearizing the Navier-Stokes equations about a steady-state flow and controlling the discrepancy between the actual flow and the steady-state flow (cf. \cite{burns2003ndp}), discretizing the associated linear state space model and then setting up the discrete flow control problem.

The fluid flow about a rotating circular cylinder can be described using the Navier-Stokes equations
\begin{eqnarray*}
  \frac{\partial\mathbf{v}}{\partial t}
  +\mathbf{v} \cdot\nabla\mathbf{v}  &=&-\nabla{p}
  +\frac{1}{Re} \nabla\cdot\tau(\mathbf{v})+ {\cal B} u,\\
  \nabla \cdot \mathbf{v}&=&0,
\end{eqnarray*}
where $\mathbf{v}$ is the fluid velocity, $p$ is the pressure, $\tau(\mathbf{v})=\nabla \mathbf{v} + \nabla \mathbf{v}^T$ is the viscous stress tensor, and 
${\cal B}u$ is the prescription of the Dirichlet boundary conditions
on the cylinder surface with $u(t)$ representing the instantaneous tangential velocity component.  Our strategy is to linearize these equations about a
desired flow profile, then use this linearized model to regulate the discrepancy between the actual flow and the desired flow.  For this study, we selected the steady-state solution at $Re=60$ created from a uniform free-stream velocity profile.  Although the steady-state solution $(\bar{\mathbf{v}},\bar{p})$ exists at this low Reynolds number, it is an unstable equilibrium solution of the Navier-Stokes equations, solving
\begin{eqnarray*}
  \bar{\mathbf{v}}\cdot\nabla\bar{\mathbf{v}} &=& -\nabla \bar{p} + \frac{1}{Re} \nabla\cdot\tau(\bar{\mathbf{v}}), \\
  \nabla \cdot \bar{\mathbf{v}}&=& 0,
\end{eqnarray*}
computed with uniform inflow velocity $\bar{\mathbf{v}}=(1,0)$ and zero velocity on the cylinder surface $\bar{\mathbf{v}}=\bar{u}\hat{t}=(0,0)$.  For computational purposes, we consider the finite flow domain $\Omega$ consisting of the unit diameter cylinder centered at the origin embedded in a rectangular flow domain $(-7,15)\times(-7,7)$.  For boundary conditions, we specify the uniform inflow velocity at the $\xi=-7$ boundary and stress-free outflow boundary conditions on the $\eta=-7$, $\xi=7$, and $\eta=7$ edges.

If we write $\mathbf{v}=\bar{\mathbf{v}}+\mathbf{v}^\prime$ and $p = \bar{p}+p^\prime$, then the flow fluctuations $(\mathbf{v}^\prime,p^\prime)$ satisfy the equations
\begin{eqnarray}
\label{eq:lin_momentum}
  \frac{\partial\mathbf{v}^\prime}{\partial t} &=&
  -\mathbf{v}^\prime\cdot\nabla\bar{\mathbf{v}} 
  -\bar{\mathbf{v}}\cdot\nabla\mathbf{v}^\prime 
  -\nabla p^\prime
  +\frac{1}{Re} \nabla\cdot\tau(\mathbf{v}^\prime) + {\cal B}u + {\bf F}(\mathbf{v}^\prime) \\ 
\label{eq:lin_cont}
  0 &=& \nabla\cdot\mathbf{v}^\prime,
\end{eqnarray}
where ${\bf F}(\mathbf{v}^\prime)$ satisfies $\| {\bf F}(\mathbf{v}^\prime) \| = O(\|\mathbf{v}^\prime\|^2)$.  The velocity fluctuation satisfies homogeneous boundary conditions at $\xi=-7$ and stress-free boundary conditions on the remaining exterior boundaries.  Using ${\cal B}u(\cdot)$ to drive $\mathbf{v}^\prime \rightarrow \mathbf{0}$ is equivalent to using ${\cal B}u(\cdot)$ to control the flow $(\mathbf{v},p)$ to $(\bar{\mathbf{v}},\bar{p})$.  We seek to achieve this using linear control theory.  Note that ignoring the nonlinear term ${\bf F}(\mathbf{v}^\prime)$ in (\ref{eq:lin_momentum}) produces the Oseen equations and often arises when developing linear feedback flow control strategies for the Navier-Stokes equations, cf. \cite{burns2003ndp}.  

At this point, we follow the standard strategy for calculating the linear feedback control laws for this problem (known as the {\em reduce-then-control} approach).  We first develop a suitable discretization for equations (\ref{eq:lin_momentum})--(\ref{eq:lin_cont}) which results in a large system of DAEs and formulate the associated linear control problem for this approximate model.  The solution to the resulting control problem is challenging and typically requires the use of suitable model reduction methods.  The presentation of a new model reduction strategy for this class of problems will then be provided in Section \ref{sec:int_dae}.

\subsection{\label{sec:lindae}Finite Element Discretization and the DAE Control Problem}
We use a standard Taylor-Hood (P2-P1) finite element pair to find approximations to both $(\bar{\mathbf{v}},\bar{p})$ and $(\mathbf{v}^\prime,p^\prime)$, cf. \cite{gunzburger1989fem}.  The nodal values of the fluctuating velocity components are denoted by $\mathbf{x}_1(t)$ while those for the pressure are denoted by $\mathbf{x}_2(t)$.  We considered several choices for the controlled output variable $y$, but for the computations below, we define 
\begin{equation}
\label{eq:C=avg_vort}
  y_{[2i-1,2i]}(t) = \frac{1}{|\Omega_i|}\int_{\Omega_i}  \mathbf{v}^\prime(t,\xi) \ d\xi \ d\eta, \qquad t>0,
\end{equation}
for six different patches downstream of the cylinder, located at $\Omega_1=[1,2.5]\times[0,2]$, $\Omega_2=[2.5,4]\times[0,2]$,
$\Omega_3=[4,5.5]\times[0,2]$, and three more reflected about the $\xi$ axis.  
For each patch, we recover two components of the average fluctuating velocity.  This is discretized as 
\begin{displaymath}
  \mathbf{y}(t) = \mathbf{C}_1 \mathbf{x}_1(t) \qquad (\mbox{and generally }
  \mathbf{y}(t) = \mathbf{C} \mathbf{x} + \mathbf{D} \mathbf{u}),
\end{displaymath}
and leads to $p=12$ controlled output variables.

We now describe the flow control problem as well as the discretized version, the DAE control problem.  The ultimate objective is to minimize the average 
fluctuation of the velocity from the smooth steady-state flow by optimally prescribing the rotational motion of the cylinder.  For well-posedness, we place a penalty on activating the control.  Thus, the control problem is
\begin{displaymath}
  \min_u \int_0^\infty \left\{ \mathbf{y}^T(t)\mathbf{y}(t) + u^T(t)R u(t) \right\} \ dt,
\end{displaymath}
where $R>0$ is a preselected constant (taken as 10 in this study), and 
subject to the constraint that the flow satisfies (\ref{eq:lin_momentum})--(\ref{eq:lin_cont}) from some initial perturbed flow state.

Upon discretization, 
the problem becomes: Find a control $\mathbf{u}(\cdot)$ that solves
\begin{equation}
\label{eq:regulator}
   \min_{\mathbf{u}}  \int_0^\infty \left\{ \mathbf{x}_1^T(t) \mathbf{C}_1^T \mathbf{C}_1  \mathbf{x}_1(t)
      + \mathbf{u}^T(t) \mathbf{R} \mathbf{u}(t) \right\}  \ dt,
\end{equation}
subject to
\begin{equation} \label{eq:dae}
      \left[ \begin{array}{cc} \mathbf{E}_{11} &  \mathbf{0}  \\ \mathbf{0} &  \mathbf{0} 
      \end{array} \right] 
      \left[ \begin{array}{c}  \dot{\mathbf{x}}_1(t)  \\  \dot{\mathbf{x}}_2(t) 
      \end{array} \right] =
      \left[ \begin{array}{cc} \mathbf{A}_{11} &  \mathbf{A}_{21}^T \\ \mathbf{A}_{21} &  \mathbf{0}
      \end{array} \right] 
      \left[ \begin{array}{c} \mathbf{x}_1(t) \\ \mathbf{x}_2(t) 
      \end{array} \right] +
      \left[ \begin{array}{c} \mathbf{B}_1 \\ \mathbf{0}
      \end{array} \right]
      \mathbf{u}(t),
\end{equation}
where $\mathbf{E}_{11} \in \bbbr^{n_1\times n_1} $ is the mass matrix for the velocity fluctuation variables and has full rank.  The matrix $\mathbf{A}_{11} \in \bbbr^{n_1\times n_1}$, $\mathbf{A}_{21}\in \bbbr^{n_2\times n_1}$, and $\mathbf{B}_1 \in \bbbr^{n_1\times m}$ (note that we use the Dirichlet map, cf. \cite{lasiecka1980uta}, and $m=1$ for this problem). Since we consider stress-free outflow conditions, 
$\bfA_{21}$  has full rank and
$\bfA_{21}\bfE_{11}^{-1}\bfA_{21}^T$ is nonsingular.  Additionally, since the tangential velocity control doesn't add mass to the domain, the term $\mathbf{B}_2$ does not appear above.

To capture the von K\'arm\'an vortex street, as well as to resolve the influence of cylinder rotations on the flow for this modest Reynolds number of 60, we use a mesh with about 5,400 elements.  However, in most flow control problems, typical dimensions of $n_1$ and $n_2$ prohibit the straight-forward application of linear control methods to the problem above.  Therefore, we investigate the use of interpolatory model reduction methods to create modest size problems from which we can develop suitable feedback control laws.

\section{Interpolatory Projections}
In this section, we describe the details of the interpolatory model reduction methodology we employ. We will explain the interpolation techniques for both the general DAE framework and the index-$2$ Oseen model arising in our application as explained in Section 
\ref{sec:flow}.

\subsection{Interpolatory Model Reduction of DAEs}
Consider the following system of differential algebraic equations (DAEs) given in  the state-space form:
\begin{equation} \label{dae_fom}
\arraycolsep 2pt
    \begin{array}{rcl}
      \bfE\, \dot{\bfx}(t) & = & \bfA\bfx(t)+\bfB\bfu(t),  \\
       \bfy(t) & = & \bfC\bfx(t) + \bfD \bfu(t),
    \end{array}  
\end{equation}
where  $\bfx(t) \in \IR^n$ represent the internal variables, $\bfu(t) \in \IR^m$ are the inputs (excitation) and $\bfy(t) \in \IR^p$ are
the outputs (observation) of the underlying dynamical system.  
In (\ref{dae_fom}), $\bfE\in \IR^{n\times n}$ is a {\it singular} matrix, thus  leading to a DAE system,
$\bfA\in \IR^{n\times n}$,  $\bfB\in \IR^{n\times m}$,  $\bfC\in  \IR^{p\times n}$, 
and $\bfD\in \IR^{p\times m}$. 
The model reduction framework for linear dynamical systems, especially for the  interpolatory methods, is best understood in the frequency domain. Towards this goal, let $\widehat{\bfu}(s)$ and $\widehat{\bfy}(s)$ denote the Laplace transforms of 
$\bfu(t)$ and $\bfy(t)$, respectively, and   take the Laplace transformation of  (\ref{dae_fom}) 
to obtain 
\begin{equation}  \label{eq:tf}
\widehat{\bfy}(s) = \bfG(s) \widehat{\bfu}(s), \mbox{~where~} \bfG(s) = \bfC(s\bfE - \bfA)^{-1}\bfB + \bfD.
\end{equation}
In (\ref{eq:tf}), $\bfG(s)$ is called the transfer function of (\ref{dae_fom}).
We will denote both the underlying dynamical system and 
its transfer function by $\bfG$.

In this setting of model reduction, the goal is to construct a reduced model of the form
 \begin{equation} \label{dae_rom}
  \arraycolsep 2pt
  \begin{array}{rcl}
       \bfEr\, \dot{\widetilde{\bfx}}(t) & = & \bfAr\widetilde{\bfx}(t)+\bfBr\bfu(t),  \\
       \widetilde{\bfy}(t) & = & \bfCr\widetilde{\bfx}(t) + \bfDr \bfu(t),
    \end{array}
\end{equation}
where $\bfEr,\bfAr\in\IR^{r\times r}$, $\bfBr\in\IR^{r\times m}$, $\bfCr\in\IR^{p\times r}$, 
and $\bfDr\in \IR^{p\times m}$ with $r \ll n$ such that 
the reduced model output 
$\widetilde{\bfy}(t)$ approximates the original output $\bfy(t)$ for a wide range of input selections 
 $\bfu(t)$ with bounded energy. 
 As for the full model, we obtain the transfer function of the reduced model by taking the Laplace transform of 
 (\ref{dae_rom}):
 \begin{equation}
 \bfGr(s)=\bfCr(s\bfEr - \bfAr)^{-1}\bfBr + \bfDr.
 \end{equation}
 Thus, in the frequency domain, we can view the model reduction problem as a rational approximation problem in which we search for a reduced   order rational function $\bfGr(s)$ to approximate the full order one $\bfG(s)$.
 
We will employ the commonly used Petrov-Galerkin projection framework to obtain the reduced model. 
We will construct  two model reduction bases  $\bfV \in \IR^{n \times r}$ and $\bfW \in \IR^{n \times r}$, 
approximate the full-order state $\bfx(t)$ by $\bfV\widetilde{\bfx}(t)$, and obtain the reduced-order model in (\ref{dae_rom}) using 
\begin{equation}  \label{red_projection}
\bfEr= \bfW^{T} \bfE \bfV, \quad\bfAr = \bfW^{T} \bfA \bfV,\quad
 \bfBr = \bfW^{T}\bfB, \quad \mbox{and } \quad\bfCr = \bfC \bfV.
\end{equation}
The feedthrough term  $\bfDr$ will be chosen appropriately to enforce matching around $s=\infty$. For the case of ordinary differential equations (ODEs) where $\bfE$ is nonsingular, the generic choice is $\bfDr = \bfD$. However, for DAEs due to the eigenvalue of the matrix pencil $\lambda \bfE - \bfA$ at infinity, special care is needed in choosing $\bfDr$.

\subsection{Model Reduction by Rational Tangential Interpolation}
\label{sec:int_dae}

In model reduction by tangential interpolation, the goal is to construct a reduced transfer function $\bfGr(s)$ that 
interpolates  $\bfG(s)$ at selected points in the complex plane along selected directions. The interpolation data consists of  the interpolation points $\{\sigma_i\}_{i=1}^r \in \IC$ together with the left  tangential directions $\{\bfsfc_i\}_{i=1}^r \in\IC^p$ and the right  tangential directions $\{\bfsfb_i\}_{i=1}^r \in\IC^m$. The usage of the terms ``left" and ``right" will be clarified once we define the interpolation problem: Given $\bfG(s)$ and the interpolation data, find  
a reduced  model $\bfGr(s)=\bfCr(s\bfEr - \bfAr)^{-1}\bfBr + \bfDr$ that satisfies, for $j=1,\ldots,r,$
\begin{equation} \label{eq:tan_int}
\begin{array}{rcll}
\bfsfc_i^T \bfG(\sigma_j) &=& \bfsfc_i^T \bfGr(\sigma_j),   \\[2mm]
\bfG(\sigma_j) \bfsfb_j &=& \bfGr(\sigma_j)\bfsfb_j,  & \mbox{and}\\[2mm]
\bfsfc_i^T\bfG'(\sigma_j) \bfsfb_j &=& \bfsfc_i^T\bfGr'(\sigma_j)\bfsfb_j. 
 \end{array}
\end{equation}
In other words, we require the reduced rational function  $\bfGr(s)$ (the reduced model) to be a bitangential Hermite interpolant the original rational function $\bfG(s)$ (the full model). One might require interpolating  higher-order derivatives of  $\bfG(s)$ as well. Moreover, one might also choose different sets of interpolation points (i.e. the right and left interpolation points) along with the left and right tangential direction vectors. For brevity of the paper, we will   consider only up to Hermite interpolation  and choose one set of interpolation points. For the details of the general case, we refer the reader to 
\cite{Ant2010imr,GSW2013,gallivan2005mrm}.

The fundamental difference between model reduction of DAEs and ODEs is that due to the eigenvalue at infinity, the transfer function of a DAE system might contain a polynomial part.  The reduced transfer function is required to exactly match the polynomial part of $\bfG(s)$; otherwise the error around $s=\infty$ can grow unbounded leading to unbounded model reduction error. Therefore, model reduction methods for DAEs aims to enforce matching of polynomial parts; see, e.g., \cite{stykel2004gramian,mehrmann2005btm,heinkenschloss2008btm,benner2006partial,GSW2013,wil2} and the references therein.

Towards this goal, let $\bfG(s)$ be additively decomposed as 
\begin{equation}
\bfG(s) = \bfGsp(s) + \bfGip(s),
\label{eq:GspP}
\end{equation}
where  $\bfGsp(s)$ is the strictly proper rational part, i.e., $\lim_{s\to \infty} \bfG(s) = 0$ and $\bfGip(s)$ is the  polynomial part of $\bfG(s)$.   We will require that the reduced transfer function $\bfGr(s)$ have exactly the same polynomial part as $\bfG(s)$, i.e.,
$$
\bfGr(s) = \bfGrsp(s) + \bfGrip(s), \mbox{~~~where~~~}\bfGrip(s) =  \bfGip(s),
$$
and $\bfGrsp(s)$ is the strictly proper rational part.
This will guarantee that  the error transfer function does not contain a polynomial part and is simply
given by 
$
\bfG_{\rm err}(s) = \bfG(s) - \bfGr(s)  =  \bfGsp(s) - \bfGrsp(s).
$
For model reduction by tangential interpolation,  \cite{GSW2013,wyatt2012issues} showed how to  construct the model reduction bases $\bfV$ and $\bfW$ so that the reduced-model of (\ref{red_projection}) satisfies the interpolation conditions (\ref{eq:tan_int}) in addition to guaranteeing $\bfGrip(s) =  \bfGip(s)$. 
As expected, the left and right deflating subspaces of the pencil $\lambda \bfE - \bfA$ corresponding to the finite and infinite eigenvalues will play a fundamental role in achieving this goal.  The next result is a special case of  Theorem 3.1 in  \cite{GSW2013} simplified to Hermite interpolation.

\begin{theorem} \label{interp_dae}
Given $\bfG(s) = \bfC(s\bfE-\bfA)^{-1}\bfB + \bfD$, let $\bfP_l$ 
and $\bfP_r$  be the spectral projectors onto the left and right deflating subspaces 
of the pencil $\lambda\bfE-\bfA$ corresponding to the finite eigenvalues. Let the columns 
of $\bfW_\infty$ and $\bfV_\infty$ span the left and right deflating subspaces of 
$\lambda\bfE-\bfA$ corresponding to the eigenvalue at infinity. Let $\sigma_i \in \IC$, for $i=1,\ldots,r$ be interpolation points
such that $\sigma_i \bfE-\bfA$ and $\sigma_i\bfEr-\bfAr$ are nonsingular. Also let $\bfsfb_i \in \IC^{m}$ and $\bfsfc_i \in \IC^{p}$ be the corresponding tangential direction vectors for $i=1,\ldots,r$.
Construct $\bfV_{\!f}$ and 
$\bfW_{\!f}$ such that 
\begin{eqnarray}
(\sigma_i\bfE-\bfA)^{-1}\bfP_l\bfB \bfsfb_i &\in {\rm Im}(\bfV_{\!f}) \mbox{~~for~~} i=1,\ldots,r,\label{eq:vf} \\
\mbox{and~~}(\sigma_i\bfE-\bfA)^{-T}\bfP_r^T\bfC^T \bfsfc_i &\in {\rm Im}(\bfW_{\!f}) \mbox{~~for~~} i=1,\ldots,r.
 \label{eq:wf}
 \end{eqnarray}
 Then with the choice of $\bfW =[\,\bfW_{\!f}, \; \bfW_\infty\,]$,   
$\bfV=[\,\bfV_{\!f},\; \bfV_\infty\,]$, and  $ \bfDr = \bfD$,
the reduced-order model $\widetilde{\bfG}(s) = \bfCr(s\bfEr - \bfAr)^{-1}\bfBr + \bfDr$ obtained 
via projection as in \textup{(\ref{red_projection})} satisfies
the bitangential Hermite interpolation conditions (\ref{eq:tan_int}) as well as $ \widetilde{\bfP}(s) = \bfP(s)$.
\end{theorem}
Even though  Theorem \ref{interp_dae} resolves the tangential interpolation problem for DAEs, it comes with a numerical caveat that  it explicitly uses the spectral projectors $\bfP_r$ and $\bfP_l$ in the model reduction step. For large-scale DAEs,  computing
$\bfP_r$ and $\bfP_l$ is, at best, very costly if not infeasible.  Therefore, it is important to construct the model reduction bases without forming $\bfP_r$ and $\bfP_l$ explicitly. Fortunately, for the Stokes-type descriptor systems of index~2,  \cite{GSW2013,wyatt2012issues} recently showed
how to apply interpolatory projections without forming $\bfP_r$ and $\bfP_l$ explicitly. This is what we consider next.
\subsection{Interpolation Theorem for Stokes-type DAEs of Index 2}  \label{sec:mor_for_index2}
Recall the linearized DAE in (\ref{eq:dae}), together with the output equation, appearing as the constraint for the optimal control problem:
\begin{eqnarray} \label{DAEsys2}
      \left[ \begin{array}{cc} \mathbf{E}_{11} &  \mathbf{0}  \\ \mathbf{0} &  \mathbf{0} 
      \end{array} \right] 
      \left[ \begin{array}{c}  \dot{\mathbf{x}}_1(t)  \\  \dot{\mathbf{x}}_2(t) 
      \end{array} \right] &=&
      \left[ \begin{array}{cc} \mathbf{A}_{11} &  \mathbf{A}_{21}^T \\ \mathbf{A}_{21} &  \mathbf{0}
      \end{array} \right] 
      \left[ \begin{array}{c} \mathbf{x}_1(t) \\ \mathbf{x}_2(t) 
      \end{array} \right] +
      \left[ \begin{array}{c} \mathbf{B}_1 \\ \mathbf{B}_2
      \end{array} \right]
      \mathbf{u}(t), \\
      \bfy(t) &=& \bfC_1\bfx_1(t) +\bfC_2\bfx_2(t) + \bfD \bfu(t),
\end{eqnarray}
where $\bfE_{11}$ is nonsingular, $\mathbf{B}_2=\mathbf{0}$, $\bfA_{21}$  has full   rank and
$\bfA_{21}\bfE_{11}^{-1}\bfA_{21}^T$ is nonsingular. In this case, system (\ref{DAEsys2}) is of index~2.  The next theorem  from \cite{GSW2013} will show how to construct a reduced model for (\ref{DAEsys2}) without requiring the deflating projectors.  For details of the derivations, we refer the reader to \cite{GSW2013}. Also, for balanced-truncation based model reduction of  (\ref{DAEsys2}), see \cite{heinkenschloss2008btm}.
\begin{theorem} \label{interp_dae_index2}
Given are the full-order DAE in (\ref{DAEsys2}), 
and the interpolation points
 $\sigma_i \in \IC$  together with the 
tangential direction vectors
 $\bfsfb_i \in \IC^{m}$ and $\bfsfc_i \in \IC^{p}$ $i=1,\ldots,r$.  Let $\bfv_i$ and $\bfw_i$ solve
\begin{equation}  \label{eq:vi}
\left[ \begin{array}{cc} \sigma_i\bfE_{11}-\bfA_{11} ~&~\bfA_{21}^T\\  \bfA_{21} & \mathbf{0} \end{array} \right]
\left[ \begin{array}{cc} \bfv_i \\  \bfz \end{array} \right]= 
\left[ \begin{array}{cc} \bfB_1\bfsfb_i \\  \mathbf{0} \end{array} \right],
\end{equation} 
and
\begin{equation}  \label{eq:wi}
\left[ \begin{array}{cc} \sigma_i\bfE_{11}^T - \bfA_{11}^T~&~ \bfA_{21}^T\\  \bfA_{21} & \mathbf{0} \end{array} \right]
\left[ \begin{array}{cc} \bfw_i \\  \bfq \end{array} \right]= 
\left[ \begin{array}{cc} \bfC^T\bfsfc_i \\  \mathbf{0} \end{array} \right].
\end{equation} 
for $i=1,\ldots,r$.
Construct 
 \begin{equation} \label{eq:VW}
 \bfV =  \left[\bfv_1,\dots, \bfv_r\right], \mbox{~~and~~}
   \bfW =  \left[\bfw_1,\dots, \bfw_r\right].
\end{equation}
 Define  $$ \bfDr =\bfD - \bfC_2(\bfA_{21}\bfE_{11}^{-1}\bfA_{21}^T)^{-1}\bfA_{21}\bfE_{11}^{-1}\bfB_1.$$
Then the reduced  model 
\begin{equation} \label{eq:ind2rom}
\bfGr(s) = \bfC\bfV(s\bfW^T\bfE_{11}\bfV - \bfW^T\bfA_{11}\bfV)^{-1}\bfW^T\bfB_1 + \bfDr.
\end{equation}
 satisfies
the bitangential Hermite interpolation conditions (\ref{eq:tan_int}) as well as $ \widetilde{\bfP}(s) = \bfP(s)$.
\end{theorem}
Note that the expensive spectral projector computations  are completely avoided; the only numerical cost is the need to solve $2r$ (sparse) linear systems arising in (\ref{eq:vi}) and (\ref{eq:wi}).  

We also note that Theorem \ref{interp_dae_index2} can be directly extended to the case where the algebraic equation in  (\ref{DAEsys2}) has the form $ \mathbf{0}  =  \bfA_{21}\bfx_1(t) + \bfB_2\bfu(t)$ with $\bfB_2\neq\mathbf{0}$. The numerical cost stays the same; see \cite{GSW2013,heinkenschloss2008btm} for details.

\section{Numerical Results}
\begin{figure}
  \centerline{
    \includegraphics[width=.49\linewidth]{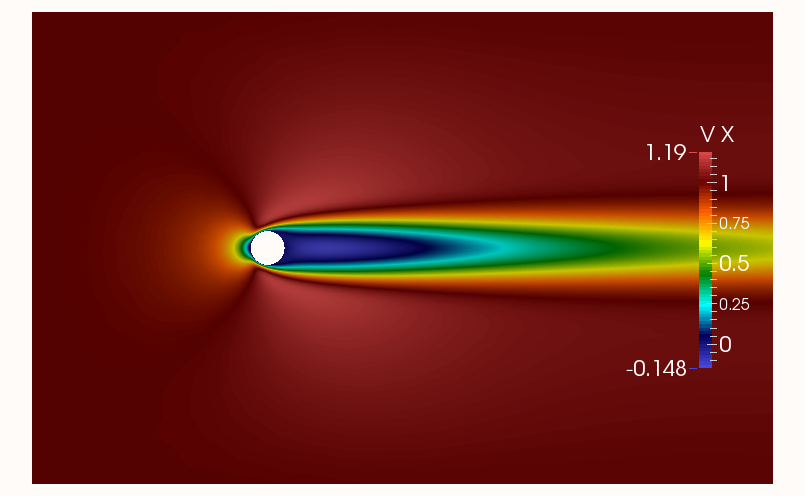}
    \includegraphics[width=.49\linewidth]{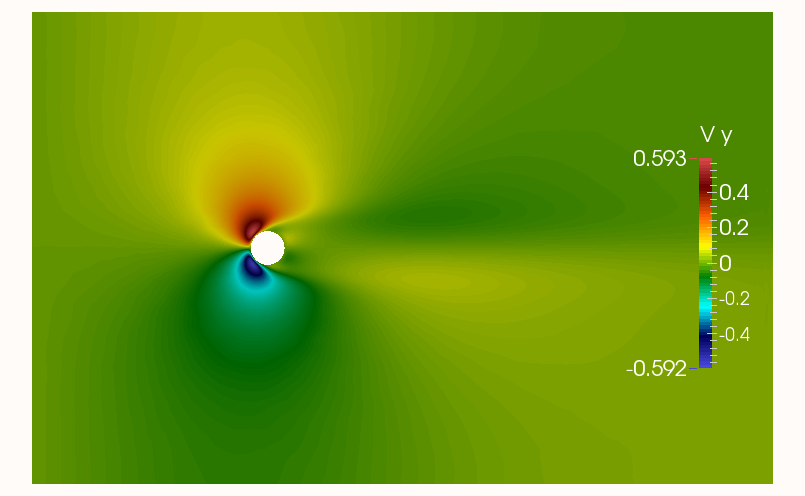}
  }
  \caption{\label{fig:v_steady}Steady-state velocity contours (horizontal-left, vertical-right)}
\end{figure}

We now apply the interpolatory model reduction algorithm described in Section \ref{sec:mor_for_index2} to the flow control problem described in Section \ref{sec:lindae}.  A relatively coarse mesh containing 5378 triangular elements was used to discretize flow solutions in the domain $\Omega = (-7,15)\times(-7,7) \backslash c$ where $c$ is the cylinder centered at the origin with unit diameter.  The steady-state flow corresponding to $Re=60$ was computed on this mesh and the resulting $(\bar{\mathbf{v}},\bar{p})$ was used to compute the discrete model for the flow fluctuations where $n_1=21,390$ and $n_2=2,777$.  Plots of $\bar{\mathbf{v}}$ components appear in Fig. \ref{fig:v_steady}.

As mentioned in Section \ref{sec:mor_for_index2}, the linearization around the steady-state solution at $Re = 60$ leads an unstable model; i.e the full-model DAE in (\ref{DAEsys2}) is unstable. There are two unstable poles are at $5.2480 \times 10^{-2} \pm \imath 7.6720 \times 10^{-1}$. We note that the model reduction framework we use does not require computing these unstable poles; we include them for comparison to the reduced model. The flexibility of the interpolatory model reduction is that even though the original model is unstable, Theorem \ref{interp_dae_index2} can  still be applied as long as the interpolation points are not chosen as one of the poles. 
We have simply chosen $33$ complex conjugate pairs (overall $66$ points) on the imaginary axis as the interpolation points. The imaginary parts of the interpolation points varied from $10^{-3}$ to $10^3$. Then, using Theorem \ref{interp_dae_index2}, we have constructed our model reduction space $\bfV  \in\IR^{21390\times 66}$. To preserve the symmetry in $\bfE_{11}$, we have simply set $\bfW = \bfV$. Thus, the reduced transfer function is a Lagrange interpolant in this case, not a Hermite interpolant. Due to the complex conjugate-pairs, construction of $\bfV$ required only $33$ sparse linear solves. Then, using a short-SVD (a relatively minor computational task due to the small number of columns in $\bfV$), we have removed the linear dependent columns from $\bfV$, and  reduced the dimension to $r=60$; thus having a final reduced model of  order $r=60$. 

\begin{figure}
\centering
 \includegraphics[scale=0.5]{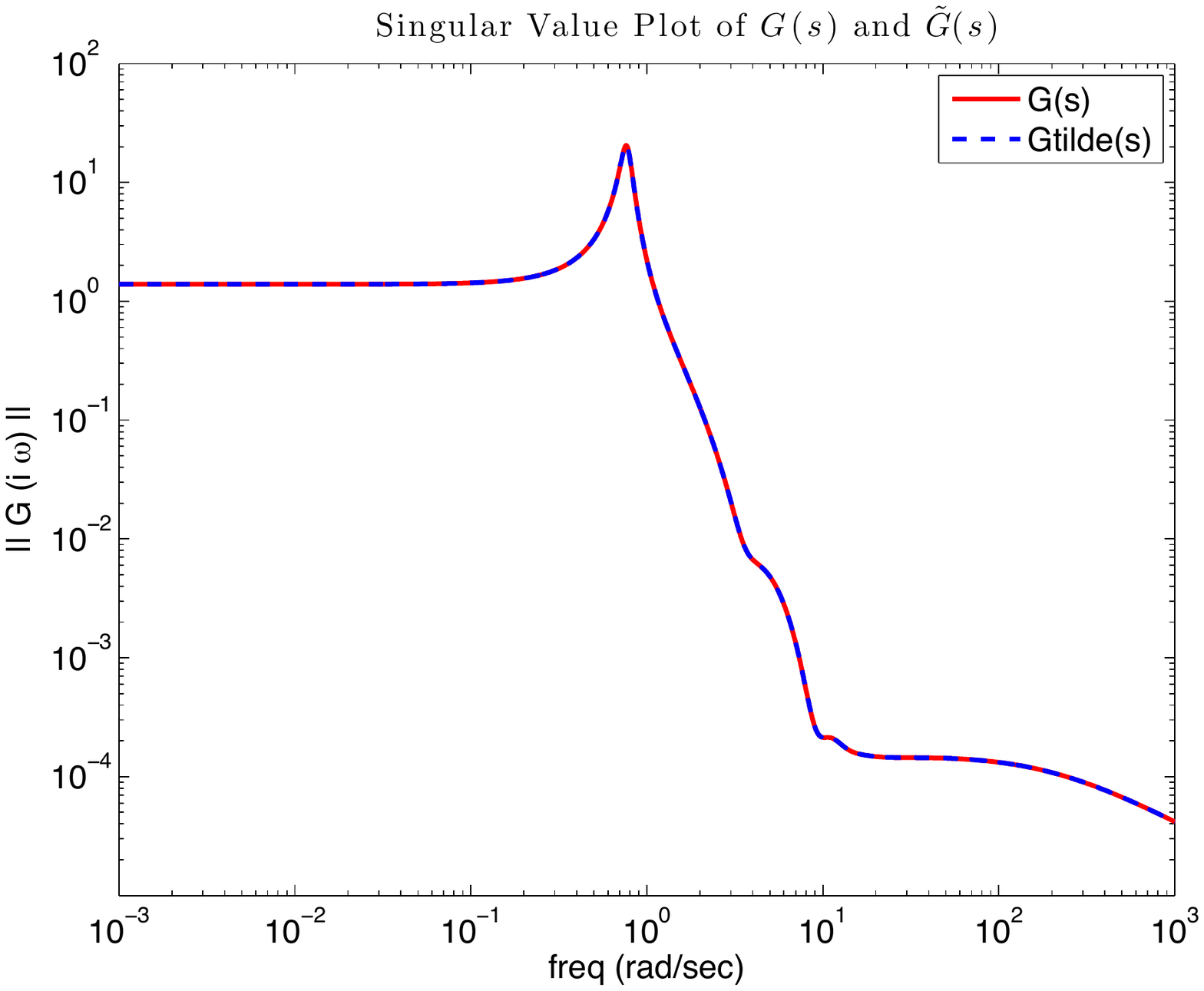} 
 \caption[ ]{The singular value  plots of the full-model $\bfG(s)$ and the reduced model $\bfGr(s)$}      
 \label{fig:Bode} 
\end{figure}

An important requirement of the reduced model in this optimal control setting is that the reduced model should capture the unstable poles of the original model so that the controller design based on the reduced model can work effectively on the full-model. As for the full-order model, our reduced model has exactly two unstable poles. The unstable poles of $\bfG(s)$ and $\bfGr(s)$ are listed below: 
$$
\begin{array}{c}
\lambda_{\rm unstable}(\bfG(s)):~~~  5.248019596820730\times 10^{-2} \pm \imath\,  7.672028760928972 \times10^{-1} \\
\lambda_{\rm unstable}(\bfGr(s)):~~~  5.248030491505502 \times10^{-2} \pm \imath\, 7.672029050490372\times10^{-1} 
\end{array}
$$
As these numbers show, the unstable poles of $\bfG(s)$ are captured to a great accuracy as desired. To further illustrate the quality of the reduced model, in Fig.~\ref{fig:Bode}, we depict the singular values plots of frequency responses of $\bfG(s)$ and $\bfGr(s)$, i.e. 
$\| \bfG(\imath \omega)\| $ and  $\| \bfGr(\imath \omega)\| $ vs $\omega \in \IR$.
As the figure shows, $\bfGr(s)$ replicates $\bfG(s)$ almost exactly.

We use the reduced-order model (\ref{eq:ind2rom}) to compute an approximate solution to the LQR problem (\ref{eq:regulator}).  Therefore, since $\widetilde{\mathbf{D}}=\mathbf{0}$ in our example, we consider the solution to the problem
\begin{equation}
\label{eq:red_lqr}
  \min_\mathbf{u} \int_0^\infty \widetilde{\mathbf{x}}_1^T(t) \widetilde{\mathbf{C}}_1^T \widetilde{\mathbf{C}}_1 \widetilde{\mathbf{x}}_1(t) + \mathbf{u}^T(t) \mathbf{R} \mathbf{u}(t)\ dt
\end{equation}
subject to $\widetilde{\mathbf{x}}_1(\cdot)$ solving (\ref{dae_rom}).  The solution can be computed by solving the algebraic Riccati equation (using the  {\tt care} function in Matlab)
\begin{displaymath}
  \widetilde{\mathbf{A}}_{11}^T \mathbf{P} \widetilde{\mathbf{E}}_{11} + \widetilde{\mathbf{E}}_{11}^T \mathbf{P} \widetilde{\mathbf{A}}_{11} - \widetilde{\mathbf{E}}_{11}^T \mathbf{P} \widetilde{\mathbf{B}}_1 \mathbf{R}^{-1} \widetilde{\mathbf{B}}_1^T \mathbf{P} \widetilde{\mathbf{E}}_{11} + \mathbf{\widetilde{C}}_1^T\mathbf{\widetilde{C}}_1 = \mathbf{0}
\end{displaymath}
for the positive definite, symmetric solution $\mathbf{P}$, then computing
\begin{displaymath}
  \widetilde{\mathbf{K}} = \mathbf{R}^{-1}\widetilde{\mathbf{B}}_1^T\mathbf{P}\widetilde{\mathbf{E}}_{11}.
\end{displaymath}
The solution to (\ref{eq:red_lqr}) is then ${\mathbf{u}} = - \widetilde{\mathbf{K}} \widetilde{\mathbf{x}}_1$.  To find the representation of the control law in the original (full-order, discrete) variables, we can use 
\begin{displaymath}
  {\mathbf{u}} = -\underbrace{\widetilde{\mathbf{K}}\mathbf{V}^T}_{\mathbf{K}} \underbrace{\mathbf{V}\widetilde{\mathbf{x}}_1}_{\approx \mathbf{x}_1}.
\end{displaymath}
Finally, we can consider the computation of $\mathbf{u}$ as an approximation of the infinite dimensional control problem where $\mathbf{v} = (v_\xi,v_\eta)$ for spatial variables $(\xi,\eta)$
\begin{displaymath}
  u(t) = -\int_\Omega h_{v_\xi}(\xi,\eta)(v_\xi(\xi,\eta,t)-\bar{v}_\xi(\xi,\eta)) + h_{v_\eta}(\xi,\eta)(v_\eta(\xi,\eta,t)-\bar{v}_\eta(\xi,\eta))\ d\xi\ d\eta.
\end{displaymath}

\begin{figure}
  \centerline{
    \includegraphics[width=.49\linewidth]{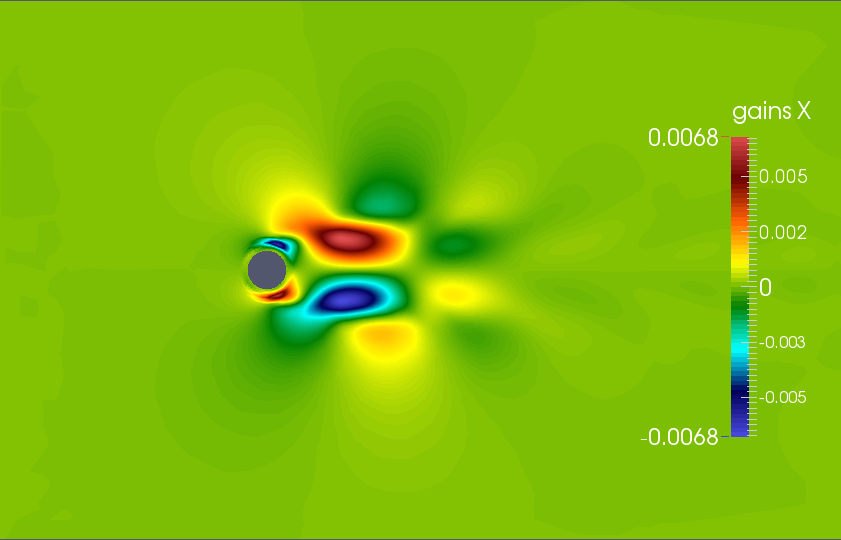}
    \includegraphics[width=.49\linewidth]{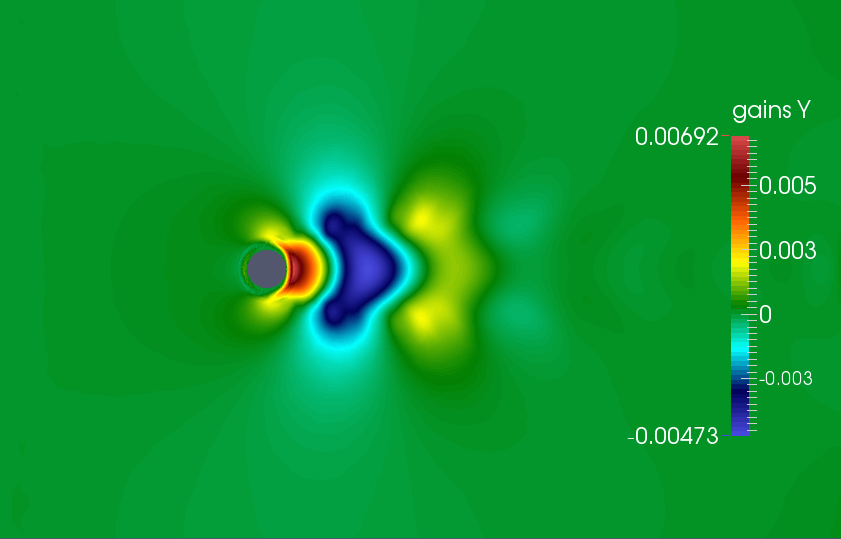}
  }
  \caption{\label{fig:f_gains}Functional gains (horizontal $h_{v_\xi}$-left, vertical $h_{v_\eta}$-right)}
\end{figure}

The finite element representations of the gains $h_{v_\xi}$ and $h_{v_\eta}$
corresponding, respectively, to the horizontal and vertical components of the velocity fluctuations are plotted in Fig.~\ref{fig:f_gains}.  These compare well to those calculated for the same problem (with a different $\mathbf{C}$ operator and slightly higher Reynolds number of 100) appearing in Fig.~4 of\cite{akhtar2010crc}.  Considering that the gains computed in this study were
computed with dramatically less computational cost emphasizes the feasibility of this approach. 

\begin{figure}
  \centerline{
    \includegraphics[width=.49\linewidth]{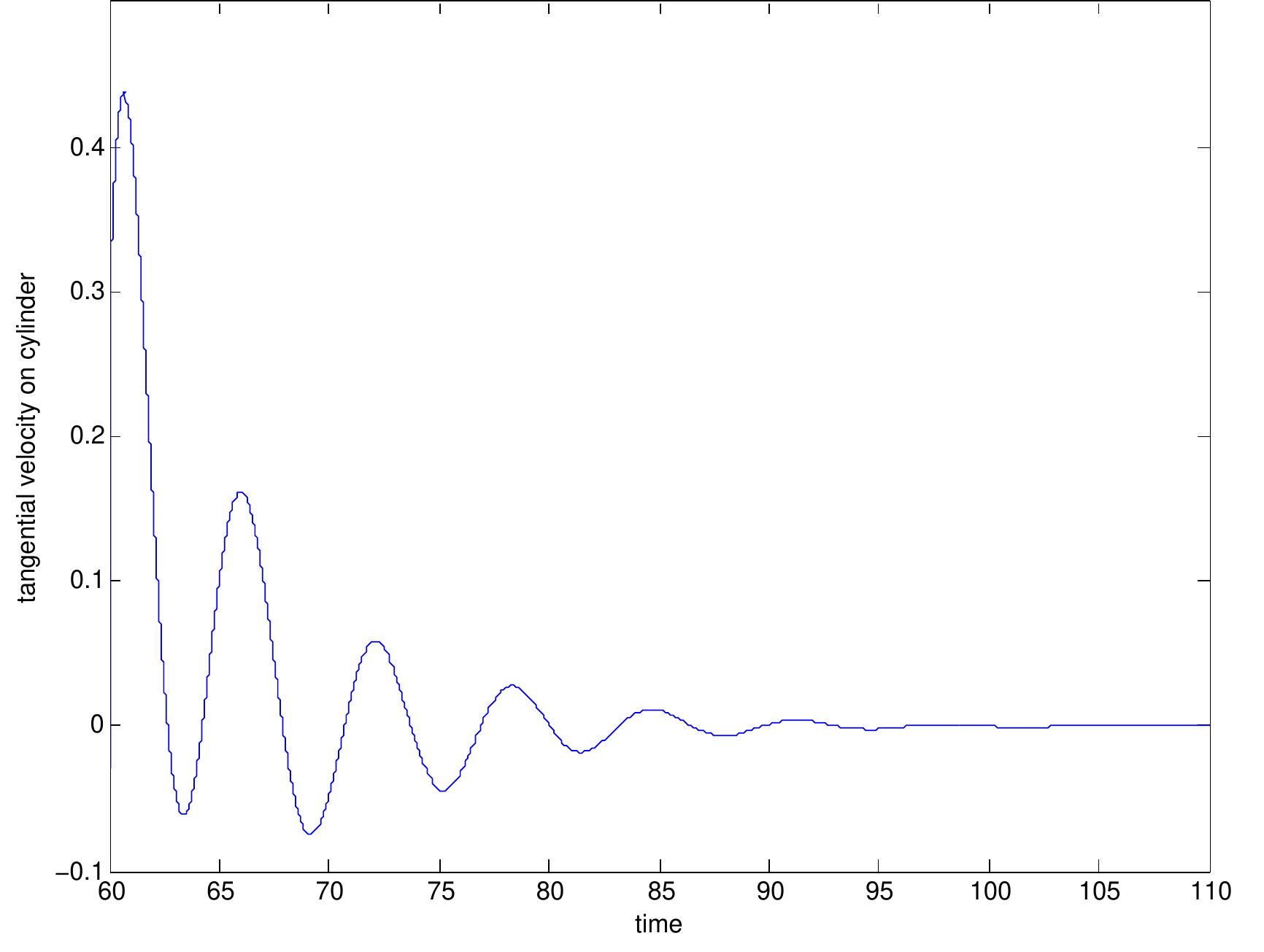}
  }
  \caption{\label{fig:control}Control Input: Tangential Cylinder Velocity}
\end{figure}

To verify that this control law stabilizes the flow, we simulated the von
K\'arm\'an vortex street for 60 seconds, then applied the full-state feedback control.  The control was able to nearly return the flow to the steady-state flow in 50 seconds of simulation.  The control input is plotted in Fig.~\ref{fig:control}.

\section{Acknowledgements}

We would like to thank Prof. Peter Benner for his valuable comments. This work was supported in part by the Air Force Office of
Scientific Research under contract FA9550-12-1-0173 and the National
Science Foundation under contracts DMS-1016450 and DMS-1217156.

\section{Conclusions and Future Work}

We have shown that with modest cost, using interpolatory model reduction we can produce accurate reduced models for index-$2$ DAEs arising from flow control problems. Using this model reduction framework within a control setting led to qualitatively similar functional gains as computed using more expensive control algorithms.  
We will apply this approach to study more complicated settings including flow control problems with higher Reynolds numbers and finer discretization. We will also investigate the practical problem of building effective state estimators and reduced-order compensators, as well as the effect of different controlled output operators on the quality of the feedback control.

\bibliographystyle{plain}

\bibliography{ActiveFlow,afc}

\end{document}